# On Least Squares Linear Regression Without Second Moment

BY RAJESHWARI MAJUMDAR

*University of Connecticut*

If $X$ and $Y$ are real valued random variables such that the first moments of $X$, $Y$, and $XY$ exist and the conditional expectation of $Y$ given $X$ is an affine function of $X$, then the intercept and slope of the conditional expectation equal the intercept and slope of the least squares linear regression function, even though $Y$ may not have a finite second moment. As a consequence, the affine in $X$ form of the conditional expectation and zero covariance imply mean independence.

## 1 Introduction

If $X$ and $Y$ are real valued random variables such that the conditional expectation of $Y$ given $X$ is an affine function of $X$, then the intercept and slope of the conditional expectation equal, respectively, the intercept and slope of the least squares linear regression function. As explained in Remark 1, when both $X$ and $Y$ have finite second moments, this equality follows from the well understood connection among conditional expectation, least squares linear regression, and the operation of projection. However, that this equality continues to hold when one only assumes that the first moments of $X$, $Y$, and $XY$ exist is the most important finding of this note [Theorem 1].

Theorem 1 evolves from the investigation of the directional hierarchy of the interdependence among the notions of independence, mean independence, and zero covariance. It is well-known that, for random variables $X$ and $Y$, independence implies mean independence and mean independence implies zero covariance, whenever the notions of mean independence and covariance make sense. Note that, the notion of covariance makes sense as long as the first moments of $X$, $Y$, and $XY$ exist; it does not require $X$ and $Y$ to have finite second moments. We review well-known counterexamples to document that the direction of this hierarchy cannot be reversed in general. Theorem 1, above and beyond establishing that an affine in $X$ form of the conditional expectation of $Y$ given $X$ implies its equality with the least squares linear regression function, also leads to the conclusion that mean independence is necessary and sufficient for zero covariance when the conditional expectation is affine in $X$.

Remark 2 examines the feasibility of obtaining the result of Theorem 1 by using the projection operator interpretation of conditional expectation and least squares linear regression elucidated in Remark 1 and the technique of extending operators to the closure of their domains, and concludes in the negative. Remark 3 explains how the relaxation of the assumption of $Y$ having a finite second moment is non-trivial.

The following notational conventions are used throughout the note. Equality (or inequality) involving measurable functions defined on a probability space, unless otherwise indicated, indicates that the relation holds almost surely. Sets are always identified with their indicator functions and $\mathfrak{K}$ denotes the universal null set. The normal distribution with mean $\mu$ and variance $\sigma^2$ is denoted by $\mathcal{N}(\mu, \sigma^2)$.

In what follows, we repeatedly use the *averaging property*, the *pull-out property*, and the *chain rule* of conditional expectation [Kallenberg (2002, page 105)].

## 2 Results

Let $X$ and $Y$ be real valued random variables on a probability space $(\Omega, \mathcal{F}, P)$. Let E denote the expectation induced by $P$ and $\mathrm{E}^{\mathfrak{C}}$ the conditional expectation given a sub $\sigma$-algebra $\mathfrak{C}$ of $\mathcal{F}$. If $\mathfrak{C} = \sigma(U)$ for some random variable $U$ defined on $(\Omega, \mathcal{F}, P)$, we write $\mathrm{E}^U$ in place of $\mathrm{E}^{\mathfrak{C}}$. Let $\mathcal{L}_1$ (respectively, $\mathcal{L}_2$) denote the Banach (respectively, Hilbert) space of integrable (respectively, square integrable) real valued functions on $(\Omega, \mathcal{F}, P)$.

If $X, Y, XY \in \mathcal{L}_1$, and $X$ and $Y$ are independent, then

$$\mathrm{Cov}(X, Y) = 0. \tag{1}$$

As is well-known, the reverse implication is not true in general; see Example 1.

**Example 1.** Let $X \sim \mathcal{N}(0, 1)$ be independent of the Rademacher random variable $W$, defined by $P(W = -1) = P(W = 1) = 1/2$. Let $Y = WX$; then $Y \sim \mathcal{N}(0, 1)$ and (1) holds. However, $X$ and $Y$ are not independent; because if they were, then $X + Y$ would have the $\mathcal{N}(0, 2)$ distribution, implying $P(X + Y = 0) = 0$, whereas in actuality $2P(X + Y = 0) = 1$. //

**Definition 1.** If $Y \in \mathcal{L}_1$ and

$$\mathrm{E}^X(Y) = \mathrm{E}(Y), \tag{2}$$

$Y$ is said to be *mean independent* of $X$. Similarly, if $X \in \mathcal{L}_1$ and

$$\mathrm{E}^Y(X) = \mathrm{E}(X), \tag{3}$$

$X$ is said to be *mean independent* of $Y$. Example 2 shows that $Y$ can be mean independent of $X$ without $X$ being mean independent of $Y$.

**Example 2.** Consider the equi-probable discrete sample space $\Omega = \{-1, 0, 1\}$ and define random variables $X$ and $Y$ on $\Omega$ as $X(\omega) = [\omega = 0]$ and $Y(\omega) = \omega$. Since $\sigma(X) = \{\mathfrak{K}, \{0\}, \{-1, 1\}, \Omega\}$ and $\mathrm{E}(YA)$ is trivially equal to 0 for all $A \in \sigma(X)$, we obtain $\mathrm{E}^X(Y) = 0 = \mathrm{E}(Y)$. However, since $X(\omega) = [Y(\omega) = 0]$, $X$ is $\sigma(Y)$ measurable, and consequently $\mathrm{E}^Y(X) = X$, whereas $\mathrm{E}(X) = 1/3$. //

It follows from the definition of independence and conditional expectation [Dudley (1989, page 264)] that for $X$ and $Y$ independent, (2) holds if $Y \in \mathcal{L}_1$, whereas (3) holds if $X \in \mathcal{L}_1$. The asymmetric nature of the notion of mean independence established in Example 2 shows that (2) (or, for that matter, (3)) does not imply independence of $X$ and $Y$. Example 3 extends Example 1 to show that even (2) and (3) combined does not necessarily imply independence of $X$ and $Y$.

**Example 3.** Let $X$, $W$, and $Y$ be as in Example 1. Since $X$ and $W$ are independent and $-X \sim \mathcal{N}(0,1)$, by Corollary 7.1.2 of Chow and Teicher (1988),

$$\mathrm{E}(X[Y \in B]) = \frac{1}{2}\mathrm{E}(X[X \in B]) + \frac{1}{2}\mathrm{E}(X[-X \in B]) = 0,$$

implying that $\mathrm{E}^Y(X) = 0 = \mathrm{E}(X)$. Clearly, by the pull-out property,

$$\mathrm{E}^X(Y) = \mathrm{E}^X(WX) = X\mathrm{E}^X(W) = X\mathrm{E}(W) = 0 = \mathrm{E}(Y).$$

However, as observed in Example 1, $X$ and $Y$ are not independent. //

As observed in the introduction, that mean independence implies zero covariance is well-known. Example 4 shows that zero covariance, that is (1), does not necessarily imply mean independence of (2).

**Example 4.** Let $X$ be uniformly distributed over the interval $(-1, 1)$ and $Y = X^2$. Since $\mathrm{E}(X) = 0 = \mathrm{E}(X^3) = \mathrm{E}(XY)$, (1) holds. Since $\mathrm{E}(Y) = \mathrm{E}(X^2) = \mathrm{Var}(X) = 1/3$ and $\mathrm{E}^X(Y) = \mathrm{E}^X(X^2) = X^2$, (2) does not hold. //

Theorem 1 shows that if $Y$ is mean independent of $X$ (as in (2)), then (1) holds; it also characterizes the setup wherein the reverse implication holds.

**Theorem 1.** Assume that $X, Y, XY \in \mathcal{L}_1$. Then (2) implies (1) and

$$\mathrm{E}^X(Y) = \alpha + \beta X \tag{4}$$

for some $\alpha, \beta \in \mathfrak{R}$. Conversely, (4) implies

$$\alpha = \mathrm{E}(Y) - \beta \mathrm{E}(X) \tag{5}$$

and

$$\beta = \frac{\mathrm{Cov}(X,Y)}{\mathrm{Var}(X)}[\mathrm{Var}(X) > 0]; \tag{6}$$

further, (6), in conjunction with (1), implies (2).

Clearly, if $X$ is mean independent of $Y$ as in (3), then (1) holds as well. Also, going back to Example 4, we now know why (1) does not imply (2) in that context; since $\mathrm{E}^X(Y) = X^2$, (4) does not hold, and by Theorem 1, (2) cannot hold.

Proof of Theorem 1. By the averaging and pull-out properties,
$$\mathrm{E}(XY) = \mathrm{E}\big(\mathrm{E}^X(XY)\big) = \mathrm{E}\big(X\mathrm{E}^X(Y)\big). \tag{7}$$

Also, since $X \in \mathcal{L}_1$, $\mathrm{Var}(X)$ is well defined, though it may be $\infty$.

Assume (2) holds. Then (1) follows from (7); also, (4) holds with $\alpha = \mathrm{E}(Y)$ and $\beta = 0$.

Conversely, assume (4) holds. Taking expectations of both sides,
$$\mathrm{E}(Y) = \alpha + \beta \mathrm{E}(X), \tag{8}$$
whence (5) follows. Thus it remains to show (6). Note that once we show (6), (1) implies $\beta = 0$, which, via (4) and (8), implies (2).

Now note that (8) implies
$$\mathrm{E}(X)\mathrm{E}(Y) = \alpha\mathrm{E}(X) + \beta(\mathrm{E}(X))^2. \tag{9}$$

By the averaging and pull-out properties, $\mathrm{E}(|XY|) = \mathrm{E}(|X|\mathrm{E}^X(|Y|))$. Since $|X\mathrm{E}^X(Y)| \leq |X|\mathrm{E}^X(|Y|)$, we obtain by (4) that $X(\alpha + \beta X) = X\mathrm{E}^X(Y) \in \mathcal{L}_1$. Also, since $X \in \mathcal{L}_1$, $\alpha X \in \mathcal{L}_1$ for every $\alpha \in \Re$; consequently,
$$\beta X^2 \in \mathcal{L}_1. \tag{10}$$

If $X \notin \mathcal{L}_2$, that is, if $\mathrm{Var}(X) = \infty$, then the right-hand side of (6) is 0 by definition and the left-hand side of (6) is 0 by (10). If $X \in \mathcal{L}_2$, that is, if $\mathrm{Var}(X) < \infty$, by (7) and (4),
$$\mathrm{E}(XY) = \mathrm{E}(X(\alpha + \beta X)) = \alpha\mathrm{E}(X) + \beta\mathrm{E}(X^2). \tag{11}$$

Subtracting (9) from (11),
$$\mathrm{Cov}(X, Y) = \beta\mathrm{Var}(X).$$

Multiplying both sides by $k$, where
$$k = \frac{[\mathrm{Var}(X) > 0]}{\mathrm{Var}(X)} = \begin{cases} 0 & \text{if } \mathrm{Var}(X) = 0 \\ \frac{1}{\mathrm{Var}(X)} & \text{if } \mathrm{Var}(X) > 0, \end{cases}$$

we obtain
$$\frac{\mathrm{Cov}(X,Y)}{\mathrm{Var}(X)}[\mathrm{Var}(X) > 0] = \beta[\mathrm{Var}(X) > 0].$$

Now note that (6) will follow once we show that $\beta[\mathrm{Var}(X) = 0] = 0$. Since, by the definition of variance, $X[\mathrm{Var}(X) = 0] = \mathrm{E}(X)[\mathrm{Var}(X) = 0]$, by (8),

$$\alpha[\text{Var}(X) = 0] + \beta X[\text{Var}(X) = 0] = \text{E}(Y)[\text{Var}(X) = 0];$$

in other words, an affine function of $X$ with slope $\beta[\text{Var}(X) = 0]$ is identically equal to a constant, showing that the slope is equal to 0, thereby completing the proof. $\square$

**Remark 1.** Is there an element of surprise in the conclusion of Theorem 1 which asserts that (4) implies (5) and (6)? The answer is no when $X, Y \in \mathcal{L}_2$, since in that case it follows from the well-understood connection (outlined below) between conditional expectation, least squares linear regression, and the operation of projection.

For a fixed real valued measurable function $X$ on $(\Omega, \mathcal{F}, P)$, let $\mathcal{L}_2(X)$ denote the Hilbert space of real valued measurable functions $f$ on $(\Omega, \sigma(X), P)$ such that $\text{E}(f^2) < \infty$. Clearly, $\mathcal{L}_2(X) \subset \mathcal{L}_2$. Let $\|\cdot\|_2$ denote the $\mathcal{L}_2$ norm on $\mathcal{L}_2$, and by inheritance, on $\mathcal{L}_2(X)$. Define

$$\mathcal{M}_2 = \{f \in \mathcal{L}_2 : \text{for some } g \in \mathcal{L}_2(X), f = g \text{ outside of a } P\text{-null set in } \mathcal{F}\}.$$

Clearly, $\mathcal{M}_2$ is a subspace of $\mathcal{L}_2$ that contains $\mathcal{L}_2(X)$. If $\{f_n : n \geq 1\}$ is a sequence in $\mathcal{M}_2$ that converges to $f \in \mathcal{L}_2$, then, since every convergent sequence is Cauchy, $\|f_n - f_m\|_2 \to 0$ as $m, n \to \infty$. By definition of $\mathcal{M}_2$, there exists a sequence $\{g_n : n \geq 1\}$ in $\mathcal{L}_2(X)$ such that $\|f_n - f_m\|_2 = \|g_n - g_m\|_2$, implying that $\{g_n : n \geq 1\}$ is a Cauchy sequence in $\mathcal{L}_2(X)$. Since $\mathcal{L}_2(X)$ is complete, there exists $g \in \mathcal{L}_2(X)$ such that $\|g_n - g\|_2 \to 0$ as $n \to \infty$. Since $\|f_n - g_n\|_2 = 0$ for every $n \geq 1$, by the triangle inequality in $\mathcal{L}_2$, $\|f - g\|_2 = 0$, showing that $f \in \mathcal{M}_2$, that is, $\mathcal{M}_2$ is closed in $\mathcal{L}_2$.

Since $\left(\text{E}^X(Y)\right)^2 \leq \text{E}^X(Y^2)$ by the Conditional Jensen's Inequality [Dudley (1989), Theorem 10.2.7)], $\text{E}^X(Y) \in \mathcal{L}_2(X)$ by the averaging property. Let $T$ denote the map from $\mathcal{L}_2$ to $\mathcal{L}_2(X)$ that takes $Y \in \mathcal{L}_2$ to $\text{E}^X(Y) \in \mathcal{L}_2(X)$. Note that, $T$ depends on $X$, but we suppress that dependence for notational convenience. Clearly, $T$ is linear. Since $\left(\text{E}^X(Y)\right)^2 \leq \text{E}^X(Y^2)$, $\|TY\|_2 \leq \|Y\|_2$ by the averaging property. Thus, $T$ is a contraction operator. For any $Y \in \mathcal{L}_2$ and $g \in \mathcal{L}_2(X)$, using the pull-out property,

$$\langle Y - TY, TY - g \rangle = 0,$$

where $\langle \cdot, \cdot \rangle$ denotes the inner product in $\mathcal{L}_2$, and consequently,

$$\|Y - g\|_2^2 = \|Y - TY\|_2^2 + \|TY - g\|_2^2,$$

showing that $TY$ is the unique minimizer of $\|Y - g\|_2$ as $g$ varies over $\mathcal{L}_2(X)$. Recall that $\mathcal{M}_2$ is a closed subspace of $\mathcal{L}_2$; let the orthogonal projection from $\mathcal{L}_2$ to $\mathcal{M}_2$ be denoted by $\Pi_{\mathcal{M}_2}$. Then, for any $Y \in \mathcal{L}_2$ and $f \in \mathcal{M}_2$,

$$\|Y - f\|_2^2 = \|Y - \Pi_{\mathcal{M}_2} Y\|_2^2 + \|\Pi_{\mathcal{M}_2} Y - f\|_2^2,$$

showing that $\Pi_{\mathcal{M}_2} Y$ is the unique minimizer of $\|Y - f\|_2$ as $f$ varies over $\mathcal{M}_2$. Since

$$\{\|Y - g\|_2 : g \in \mathcal{L}_2(X)\} = \{\|Y - f\|_2 : f \in \mathcal{M}_2\},$$

$\Pi_{\mathcal{M}_2} Y$ equals $TY$ outside of a $P$-null set in $\mathcal{F}$.

Note that if the probability space $(\Omega, \sigma(X), P)$ is complete, so that the almost sure limit of a sequence of measurable functions is measurable, $\mathcal{L}_2(X)$ becomes a closed subspace of $\mathcal{L}_2$, and in our identification of the conditional expectation as a projection, we can avoid the construction involving $\mathcal{M}_2$.

Let $\mathcal{H}$ denote the two-dimensional linear space spanned by $J$ and $X$, where $J$ is the real valued function defined on $(\Omega, \mathcal{F}, P)$ that is almost surely equal to 1; since $X \in \mathcal{L}_2$, we obtain $\mathcal{H} \subset \mathcal{L}_2(X) \subset \mathcal{M}_2$. Putting the degenerate case when $X$ is a scalar multiple of $J$ aside for the time being and applying the Gram-Schmidt orthonormalization process to the basis $\{J, X\}$ of $\mathcal{H}$, we obtain that $\{J, X^*\}$, where

$$X^* = \frac{X - \langle J, X \rangle J}{\|X - \langle J, X \rangle J\|_2},$$

is an orthonormal basis of $\mathcal{H}$.

Let $\mathcal{A}$ denote the subspace of $\mathcal{L}_2$ that consists of all $Y \in \mathcal{L}_2$ such that $TY \in \mathcal{H}$, that is, (4) holds for $Y$. If $Y \in \mathcal{A}$, then $\Pi_{\mathcal{M}_2} Y = TY \in \mathcal{H} \subset \mathcal{L}_2(X) \subset \mathcal{M}_2$, and consequently,

$$TY = \Pi_{\mathcal{M}_2} Y = \Pi_{\mathcal{H}}(\Pi_{\mathcal{M}_2} Y) = \Pi_{\mathcal{H}} Y = \langle J, Y \rangle J + \langle X^*, Y \rangle X^*, \tag{12}$$

where $\Pi_{\mathcal{H}}$ denotes the orthogonal projection from $\mathcal{L}_2$ to $\mathcal{H}$. In more familiar notation, (12) asserts

$$\mathrm{E}^X(Y) = \mathrm{E}(Y) + \frac{\mathrm{Cov}(X, Y)}{\mathrm{Var}(X)}(X - \mathrm{E}(X)),$$

that is, (5) and (6) hold. Note that, if $X$ is a scalar multiple of $J$, that is, $X$ is a degenerate random variable and $\mathrm{Var}(X) = 0$, then $\mathcal{H}$ is simply the span of $J$, and we obtain that $TY = \langle J, Y \rangle J$, that is, $\mathrm{E}^X(Y) = \mathrm{E}(Y)$, leading to the conclusion of Theorem 1 for $X, Y \in \mathcal{L}_2$. //

**Remark 2.** Can the conclusion of Theorem 1, when we only have $Y \in \mathcal{L}_1 \setminus \mathcal{L}_2$, $X \in \mathcal{L}_1$, and $XY \in \mathcal{L}_1$, be obtained using the projection operator tools employed in Remark 1? The answer, as far as we can tell, is no.

The domain of the map $T$ can be expanded to $\mathcal{L}_1$, causing the range to be expanded to $\mathcal{L}_1(X)$, the Banach space of real valued measurable functions $f$ on $(\Omega, \sigma(X), P)$ such that $\mathrm{E}(|f|) < \infty$. The linearity of $T$ is not impacted by this expansion of domain. Since $|\mathrm{E}^X(Y)| \leq \mathrm{E}^X(|Y|)$ (again by the Conditional Jensen's Inequality), $\|TY\|_1 \leq \|Y\|_1$ by the averaging property, where $\|\cdot\|_1$ denotes the $\mathcal{L}_1$ norm on $\mathcal{L}_1$, and by inheritance, on $\mathcal{L}_1(X)$. Thus, $T$ remains a contraction operator. In fact, as observed by Kallenberg

(2002), $T$ on $\mathcal{L}_2$ is uniformly $\mathcal{L}_1$-continuous and its extension to a linear and continuous map on $\mathcal{L}_1$ is unique up to almost sure equivalence.

The definition of $\mathcal{A}$ can be extended to denote the subspace of $\mathcal{L}_1$ that consists of all $Y \in \mathcal{L}_1$ such that $TY \in \mathcal{H}$. Since $\mathcal{H}$ is a finite-dimensional, hence closed, subspace of $\mathcal{L}_1$, and $T$ is continuous, $\mathcal{A}$ is a closed subspace of $\mathcal{L}_1$.

Let $\mathcal{B}$ denote the subspace of $\mathcal{L}_1$ that consists of all $Y \in \mathcal{L}_1$ such that $XY \in \mathcal{L}_1$. If $X \in \mathcal{L}_2$, then, by the Cauchy-Schwartz inequality, $\mathcal{L}_2 \subset \mathcal{B}$. Theorem 1 asserts that the representation of $T$ as the orthogonal projection to $\mathcal{H}$ that is valid on $\mathcal{A} \cap \mathcal{L}_2$ can be extended to hold on $\mathcal{A} \cap \mathcal{B}$. Given that the set inclusion relationship between the closure of $\mathcal{A} \cap \mathcal{L}_2$ in $\mathcal{L}_1$ and $\mathcal{A} \cap \mathcal{B}$ is generally indeterminable (if $X$ is bounded, then $\mathcal{B} = \mathcal{L}_1$ and $\mathcal{A} = \mathcal{A} \cap \mathcal{B}$, implying that the closure of $\mathcal{A} \cap \mathcal{L}_2$ in $\mathcal{L}_1$ is contained in $\mathcal{A} \cap \mathcal{B}$, though the inclusion cannot be necessarily reversed), it seems that the result of Theorem 1, even when $X$ is bounded, can not be deduced using the technique of operator extension.

What happens if $X \in \mathcal{L}_1 \backslash \mathcal{L}_2$? Clearly, $\mathcal{H}$ is no longer a subspace of $\mathcal{L}_2$. However, since $\mathcal{L}_2$ is dense in $\mathcal{L}_1$, we can find a sequence $\{X_n : n \geq 1\}$ in $\mathcal{L}_2$ such that $X_n$ converges to $X$ in $\mathcal{L}_1$. With $\mathcal{H}_n$ denoting the span of $J$ and $X_n$, the orthogonal projection of $\mathcal{L}_2$ to $\mathcal{H}_n$ converges to the orthogonal projection of $\mathcal{L}_2$ to the span of $J$. From the proof of Theorem 1, $T$ on $\mathcal{A} \cap \mathcal{B}$ equals the "orthogonal projection," in a limiting sense, to the span of $J$. //

**Remark 3.** What does the relaxation of the structural assumption from $X, Y \in \mathcal{L}_2$ to $Y \in \mathcal{L}_1 \backslash \mathcal{L}_2$, $X \in \mathcal{L}_1$, and $XY \in \mathcal{L}_1$ entail? One can conceivably argue that the Cauchy-Schwartz inequality remains the primary tool for verifying that $XY \in \mathcal{L}_1$ when $X$ and $Y$ are dependent random variables, and as such, this relaxation of assumption is neither insightful nor useful. While that argument may have some merit, we would like to point out that if $X$ is bounded, then obviously $X \in \mathcal{L}_1$, and $Y \in \mathcal{L}_1 \backslash \mathcal{L}_2$ implies $XY \in \mathcal{L}_1$.

A classic example of (4) holding for a bounded random variable $X$ is the Bernoulli random variable, since, for any measurable function $h(X)$ of the Bernoulli random variable $X$, we have $h(X) = h(0) + (h(1) - h(0))X$, and $\mathrm{E}^X(Y)$ is a measurable function of $X$.

Working with the non-central $t$-distribution with 2 degrees of freedom, the following example presents $X, Y$ such that $X \in \mathcal{L}_1$ but is not bounded, $Y \in \mathcal{L}_1 \backslash \mathcal{L}_2$, $XY \in \mathcal{L}_1$, and (4) holds. Let $X \sim \mathcal{N}(0, 1)$ and given $X$, $W \sim \mathcal{N}(X, 1)$. Let $V \sim \chi_2^2$ be independent of $(W, X)$. Let $Y = W/\sqrt{V/2}$.

Clearly, $X \in \mathcal{L}_1$ but is not bounded.

To verify that $Y \in \mathcal{L}_1 \backslash \mathcal{L}_2$, we first obtain the marginal distribution of $W$. Using the form of the conditional density of $W$ given $X = x$ and the form of the marginal density of $X$, we obtain that the joint density of $(W, X)$ is given by

$$f(w,x) = \frac{\exp\left(-\frac{w^2}{2} + wx - x^2\right)}{2\pi};$$

since $x \mapsto \exp\left(-(x-w/2)^2\right)/\sqrt{\pi}$ represents the density of the $\mathcal{N}\left(\frac{w}{2}, \frac{1}{2}\right)$ distribution, completing the square in $x$ we obtain the marginal density of $W$ to be

$$f_W(w) = \int_\Re f(w,x)dx = \frac{\exp\left(-\frac{w^2}{4}\right)}{\sqrt{2\pi}\sqrt{2}} \int_\Re \frac{\exp\left(-\left(x-\frac{w}{2}\right)^2\right)}{\sqrt{\pi}} dx = \frac{\exp\left(-\frac{w^2}{4}\right)}{\sqrt{2\pi}\sqrt{2}},$$

showing that $W \sim \mathcal{N}(0,2)$. Since $V = 2U$, where $U$ is a $\Gamma(1)$ random variable [Fabian and Hannan (1985, Definition 3.3.6)], by Theorem 3.3.9 of Fabian and Hannan (1985),

$$\mathrm{E}\left(\frac{1}{\sqrt{V/2}}\right) = \mathrm{E}\left(U^{-\frac{1}{2}}\right) = \frac{\Gamma\left(\frac{1}{2}\right)}{\Gamma(1)} = \sqrt{\pi}; \tag{13}$$

since $W$ and $V$ are independent, $Y \in \mathcal{L}_1$. However, $\mathrm{E}(2/V) = \mathrm{E}(U^{-1}) = \infty$, showing that $Y \notin \mathcal{L}_2$.

To show that $XY \in \mathcal{L}_1$, by independence of $(W,X)$ and $V$ along with (13), it suffices to show that $XW \in \mathcal{L}_1$, which follows readily from the Cauchy-Schwartz inequality.

Finally, by Corollary 7.1.2 of Chow and Teicher (1988) and (13),

$$\mathrm{E}^{W,X}(Y) = \sqrt{\pi}W,$$

implying, by the chain rule, $\mathrm{E}^X(Y) = \sqrt{\pi}X$, that is, (4) holds. //

RAJESHWARI MAJUMDAR
rajeshwari.majumdar@uconn.edu
PO Box 47
Coventry, CT 06238